\documentclass[12pt]{amsart}
\usepackage[T1]{fontenc}
\usepackage[english]{babel}
\usepackage{mathtools}
\usepackage{amstext, amssymb}
\usepackage{mathrsfs}
\usepackage{nicematrix}
\usepackage{slashed}
\usepackage{accents}
\usepackage{enumerate}
\usepackage{xfrac}
\usepackage[all]{xy}
\usepackage{tikz}
\usetikzlibrary{arrows,decorations.pathmorphing,backgrounds,positioning,fit,petri,decorations}
\usetikzlibrary{calc,intersections,through,backgrounds,mindmap,patterns,fadings}
\usetikzlibrary{decorations.text}
\usetikzlibrary{decorations.fractals}
\usetikzlibrary{fadings}
\usetikzlibrary{shadings}
\usetikzlibrary{shadows}
\usetikzlibrary{shapes.geometric}
\usetikzlibrary{shapes.callouts}
\usetikzlibrary{shapes.misc}
\usetikzlibrary{spy}
\usetikzlibrary{topaths}
\theoremstyle{plain}
\newtheorem{theorem}{Theorem}[section]
\newtheorem{proposition}[theorem]{Proposition}

\newtheorem{lemma}[theorem]{Lemma}

\theoremstyle{definition}

\newtheorem{example}[theorem]{Example}

\numberwithin{equation}{section}

\newcommand{\f}[2]{\frac{#1}{#2}}
\renewcommand{\c}{\cdot}
\renewcommand{\t}[1]{\text{#1}}
\newcommand{\alp}{\alpha}
\newcommand{\real}{\mathbb R}

\newcommand{\nat}{\mathbb N}
\newcommand{\com}[1]{\ignorespaces}
\newcommand{\pmat}[1]{\begin{pmatrix}#1\end{pmatrix}}
\newcommand{\cp}{\mathcal P}

\newcommand{\tr}{\t{tr}}

\newcommand{\lb}{\left(}
\newcommand{\rb}{\right)}

\newcommand{\IMM}{\mathscr{M}}

\newcommand{\IN}{\mathbb{N}}

\newcommand{\IR}{\mathbb{R}}

\newcommand{\Id}{\;\mathrm{d} }

\newcommand{\nn}{\nonumber}

\def\XXint#1#2#3{{\setbox0=\hbox{$#1{#2#3}{\int}$} \vcenter{\hbox{$#2#3$}}\kern-.5\wd0}}

\newcommand{\getcolor}[1]{%
    \pgfmathsetmacro{\index}{int(#1)}%
    \ifnum\index=1 \definecolor{tempcolor}{RGB}{255,80,96}
    \else\ifnum\index=2 \definecolor{tempcolor}{RGB}{0,0,255}
    \else\ifnum\index=3 \definecolor{tempcolor}{RGB}{0,255,96}
    \else\definecolor{tempcolor}{RGB}{0, 0, 0}
    \fi\fi\fi
}

\newcommand{\addnum}[2]{%
  \the\numexpr #1 + #2\relax%
}
\newcommand{\multiplynum}[2]{%
  \the\numexpr #1 * #2\relax%
}

\newcommand{\threevertexgraph}[9]{
\def\loopv{#1}
\def\loopw{#2}
\def\loopx{#3}
\def\firstvw{#4}
\def\secondvw{#5}
\def\firstvx{#6}
\def\secondvx{#7}
\def\firstwx{#8}
\def\secondwx{#9}
\threevertexgraphcont
}
\newcommand{\threevertexgraphcont}[1]{\begin{center}
\begin{tikzpicture}
    \node[circle, draw, inner sep=2pt] (v) at (0, 0) {};
    \node[circle, draw, inner sep=2pt] (w) at (3, 0) {};
    \node[circle, draw, inner sep=2pt] (x) at (6, 0) {};
    \node[below=5pt of v] {$v$};
    \node[below=5pt of w] {$w$};
    \node[below=5pt of x] {$x$};
    
    \ifnum\loopv=0\else\getcolor{\loopv}
    \draw[-, tempcolor] (v) .. controls (-1,1.5) and (1,1.5) .. (v);
    \fi
    \ifnum\loopw=0\else\getcolor{\loopw}
    \draw[-, tempcolor] (w) .. controls (2,1.5) and (4,1.5) .. (w);
    \fi
    \ifnum\loopx=0\else\getcolor{\loopx}
    \draw[-, tempcolor] (x) .. controls (5,1.5) and (7,1.5) .. (x);
    \fi
    
    \ifnum\multiplynum{\firstvw}{\secondvw}=0 \getcolor{\addnum{\firstvw}{\secondvw}}
    \def\temp{\addnum{\firstvw}{\secondvw}}
    \ifnum\temp=0\else
    \draw[-, tempcolor] (v) -- (w);
    \fi
    \else
    \getcolor{\firstvw}
    \draw[-, tempcolor] (v) .. controls (1,0.5) and (2,0.5) .. (w);
    \getcolor{\secondvw}
    \draw[-, tempcolor] (v) .. controls (1,-0.5) and (2,-0.5) .. (w);
    \fi
    
    \ifnum\multiplynum{\firstwx}{\secondwx}=0 \getcolor{\addnum{\firstwx}{\secondwx}}
    \def\temp{\addnum{\firstwx}{\secondwx}}
    \ifnum\temp=0\else
    \draw[-, tempcolor] (w) -- (x);
    \fi
    \else
    \getcolor{\firstwx}
    \draw[-, tempcolor] (w) .. controls (4,0.5) and (5,0.5) .. (x);
    \getcolor{\secondwx}
    \draw[-, tempcolor] (w) .. controls (4,-0.5) and (5,-0.5) .. (x);
    \fi
    
    \ifnum\multiplynum{\firstvx}{\secondvx}=0 \getcolor{\addnum{\firstvx}{\secondvx}}
    \def\temp{\addnum{\firstvx}{\secondvx}}
    \ifnum\temp=0\else
    \draw[-, tempcolor] (v) .. controls (2,-1.5) and (4,-1.5) .. (x);
    \fi
    \else
    \getcolor{\firstvx}
    \draw[-, tempcolor] (v) .. controls (2,-1.5) and (4,-1.5) .. (x);
    \getcolor{\secondvx}
    \draw[-, tempcolor] (v) .. controls (2,-2) and (4,-2) .. (x);
    \fi
    
    \node at (-1,0) {$#1$:};
\end{tikzpicture}
\end{center}}

\allowdisplaybreaks

\title[]{Multigraphs and Time Ordered Isserlis-Wick formulae}

\author{Sergio Cacciatori}

\author{Batu Güneysu}

\author{Sebastian Wündsch}

\begin{document}
	\begin{abstract}Given a Gaussian process $X(s)\in\mathbb{R}^m$, $s\in [0,1]$, having a scalar covariance matrix, and a polynomial $Q$ of $m$ variables with real coefficients, we calculate the expression 
    $$
    \int_{\{0\leq s_1\leq...\leq s_n\leq 1\}}\mathbb{E}\left[\prod_{k=1}^nQ(X(s_k))\right]\;\mathrm{d}s_1\cdots\mathrm{d}s_n
    $$ 
    using labeled multigraphs with explicit symmetry factors.   
	\end{abstract}
	\maketitle

\section{Introduction}

A fundamental result in the context of Gaussian random vectors is Isserlis' Theorem \cite{isserlis} from 1918: 

\begin{theorem} If $Z=(Z_1,\dots,Z_m)\in\mathbb{R}^m$ is a zero mean Gaussian random vector, then one has
$$
\mathbb{E}\left[\prod^m_{j=1}Z_j\right]=\begin{cases}0, \text{if $m$ is odd},\\
\sum_{\text{\emph{pairings} $\mathcal{P}$ \emph{of} $\{1,\dots,m\}$}}\;
\prod_{\{i,k\}\in \mathcal{P}}\mathrm{Cov}(Z_i,Z_k),\quad\text{if $m$ is even}.
\end{cases}
$$

\end{theorem}

The combinatorics of Isserlis' Theorem has been famously rediscovered in the context of quantum field theory by Wick in 1950 \cite{wick}, which is why results of this genre are often also called Wick's Theorem or Isserlis-Wick Theorem. We refer the reader to \cite{munt} for a recent modern proof of Isserlis' Theorem and remark also that there is a more general result in this context which also works for non-Gaussian random variables: the moment-cumulants formula \cite{cumulant} (see \cite{peccati} for a modern treatment of this result).\vspace{2mm}

Assume now 
$$
X(s)=(X_1(s),\dots,X_m(s))\in\mathbb{R}^m, \quad s\in [0,1],
$$
is a zero mean Gaussian process such that the covariance matrix is of the form
$$
\mathrm{Cov}(X_i(s),X_j(t))=f(s,t)\delta_{ij},\quad s,t\in [0,1],
$$
for some jointly continuous function
$$
f:[0,1]\times [0,1]\longrightarrow \IR.
$$
Given a polynomial 
$$
Q(x)=\sum_{\alpha\in\mathbb{N}^m_0} q_\alpha x^\alpha\quad\text{of $m$ variables with real coefficients},
$$
the central mathematical expressions of this paper are expressions of the form
\begin{align}\label{opo}
   \int_{\{0\leq s_1\leq...\leq s_n\leq 1\}}\mathbb{E}\left[\prod_{k=1}^nQ(X(s_k))\right]\;\mathrm{d}s_1\cdots\mathrm{d}s_n.
\end{align}

These expressions arise naturally in quantum mechanics: For $X^{x,y;t}(s)$, $s\in [0,t]$, a Brownian bridge from $x$ to $y$ with lifetime $t$, via the Feynman-Kac formula these time ordered expressions appear explicitly in the context of heat kernels of Schrödinger operators of the form $H_Q=-\Delta/2+Q$. Namely, one has
$$
\exp(-tH_Q)(x,y)=(2\pi t)^{-m/2}\exp\left(\frac{-|x-y|^2}{2t}\right)\mathbb{E}\left[\exp\left(-\int^t_0 Q(X^{x,y;t}(s))\Id s\right)\right],
$$
which connects to (\ref{opo}), if one formally expands 
$$
\mathbb{E}\left[\exp\left(-\int^t_0 Q(X^{x,y;t}(s))\Id s\right)\right]=\sum^\infty_{n=0}  \int_{\{0\leq s_1\leq...\leq s_n\leq t\}}\mathbb{E}\left[\prod_{k=1}^nQ(X^{x,y;t}(s_k))\right]\;\mathrm{d}s_1\cdots\mathrm{d}s_n.
$$
Note also that, remarkably, the Feynman-Kac formula holds (for small times) even for certain $Q$'s like 
\begin{align}\label{qay}
Q(x)=(x,Dx)+(c,x),\quad\text{where $D\in \mathbb{R}^{m\times m}$ is symmetric and $c\in\mathbb{R}^m$,}
\end{align} 
whose Schrödinger semigroup $\exp(-tH_Q)$ is unbounded \cite{simon,broderix}.\vspace{2mm}

Our aim is to calculate (\ref{opo}) solely in terms of the covariance function $f$. We illustrate our main result with a simple example: 

\begin{example} Assume $m=1$, $n=2l$ is even, and $Q(x)=x$. Then applying the Isserlis-Wick Theorem to $\mathbb{E}\left[\prod_{k=1}^n X(s_k)\right]$ shows
\begin{align*}
  &\int_{\{0\leq s_1\leq...\leq s_n\leq 1\}}\mathbb{E}\left[\prod_{k=1}^n X(s_k)\right]\;\mathrm{d}s_1\cdots\mathrm{d}s_n=\int_{\{0\leq s_1\leq \dots \leq s_{n}\leq 1\}} \sum_{\substack{\text{pairings $P$}\\
\text{ of $\{1,\dots,n\}$} }}\prod_{\{i,j\}\in \mathcal{P}}f(s_i,s_j) \Id s_1\cdots \Id s_{n}.
\end{align*}
As for every symmetric $n\times n$ matrix $A$ one has
\begin{align}\label{haf-n}
\sum_{\substack{\text{pairings } \mathcal{P}\\
\text{of } \{1,\dots,n\}}}
\prod_{\{i,j\}\in P} A_{i,j}
=
\frac{1}{l!\,2^{l}}
\sum_{\sigma\in S_{n}}
A_{\sigma(1),\sigma(2)}\cdots A_{\sigma(n-1),\sigma(n)} ,
\end{align}
where $S_n$ denotes the permutation group of $\{1,\dots,n\}$, we can perform the deterministic calculation
\begin{align*}
&=
\frac{1}{l!\,2^{l}}
\sum_{\sigma\in S_n}
\int_{\{0\le s_1\le \dots \le s_n \le 1\}}
f(s_{\sigma(1)},s_{\sigma(2)})\cdots
f(s_{\sigma(n-1)},s_{\sigma(n)})
\, \Id s_1\cdots \Id s_n
\\
&=
\frac{1}{l!\,2^{l}}
\int_{[0,1]^n}
f(s_1,s_2)\cdots f(s_{n-1},s_n)\,
\Id s_1\cdots \Id s_n
\\
&=
\frac{1}{l!\,2^{l}}
\left(
\int_{[0,1]^2} f(s_1,s_2)\, \Id s_1 \Id s_2
\right)^{l}.
\end{align*}
\end{example}

For higher degree polynomials, one naturally has to consider pairings of multisets rather than sets in the Wick-Isserlis Theorem, in order to take several occurrences of the same random variable such as e.g. $X(s_j)^2$ into account. The ultimate combinatorics via Isserlis-Wick pairs becomes very involved, and we follow the general philosophy from perturbative quantum field theory which is to tame this combinatorics via multigraphs. Our main result is Theorem \ref{prop1gen}, which states that with some fixed set of vertices $\{v_1,\dots, v_n\}$, one has  
\begin{align}\label{klq}
&\int_{\{0\leq s_1\leq...\leq s_n\leq 1\}}\mathbb{E}\left[\prod_{k=1}^nQ(X(s_k))\right]\;\mathrm{d}s_1\cdots\mathrm{d}s_n\\
&=\f1{n!}\sum_{\substack{\alpha^{(1)}, \dots, \alpha^{(n)}\in\IN^m_0 \\\nn\frac{1}{2}\sum^n_{k=1}\sum^m_{i=1}\alpha^{(k)}_i\,\in\,\IN}}\left(\prod_{k=1}^n q_{\alpha^{(k)}}\right) \sum_{\Gamma_1}\dots\sum_{\Gamma_m}\left(\prod_{q=1}^m C(\Gamma_q)\right)\left(\prod_{\Lambda\in \mathscr{C}(\Gamma_1+\cdots+\Gamma_m)}\int_{\Lambda}  f\right),
\end{align}
where 
\begin{itemize}
    \item the $l$-th sum is over all multigraphs $\Gamma_l$ with vertex set $\{v_1,\dots, v_n\}$ and $\deg(v_j)=\alp_1^{(j)}$ for all $j\in\{1,\dots,n\}$,
\item $C(\Gamma_q)\in\mathbb{N}$ is an explicitly given combinatoric number (cf. (\ref{edya})),
\item $\mathscr{C}(\Gamma_1+\cdots+\Gamma_m)$ denotes the connected components of the multigraph $\Gamma_1+\cdots+\Gamma_m$,
\item $\int_{\Lambda}  f$ is the integral of the covariance along the multigraph $\int_{\Lambda}  f$.
\end{itemize}
As indicated above, to obtain this result, we first calculate the expectation on the LHS of (\ref{klq}) for fixed times using Isserlis-Wick and then calculate the remaining deterministic time ordered integral using multigraph labelings.\vspace{2mm}

While we are not aware of an explicit calculation of the expression (\ref{opo}) in this generality in the mathematics literature, in contrast to the usual line of arguments in perturbative quantum field theory, we use labeled multigraphs as opposed to unlabeled multigraphs (modulo automorphisms), and our symmetry factors $C(\Gamma_q)$ may be viewed as a refinement of the usual implicit symmetry factors that are given as cardinalities of automorphism groups. In fact, we demonstrate how our machinery works in Example \ref{edxaa}, by calculating (\ref{opo}) for the potential (\ref{qay}) up to order $n=3$.\vspace{4mm}

\emph{Acknowledgements:} The second named author (B.G.) would like to dedicate this work to Anton Thalmaier, who has been a constant source of support and inspiration. 

\section{Main Results}

We first introduce the basic combinatoric concepts of this paper:\vspace{1mm}

A \emph{pairing} $\mathcal{P}$ of a set $M$ with even cardinality $l$ is a set  
$$
\mathcal{P}=\{\{u_1,v_1\},\dots,\{u_{l/2},v_{l/2}\} \}
$$
given by pairwise disjoint two-element subsets of $M$. 

A \emph{multiset} is a pair $\mathscr{M}=(M,\mathsf{m})$ given by a set $M$, the so called \emph{universe of $\mathscr{M}$}, and a map $\mathsf{m}:M\to \mathbb{N}_0$, called the \emph{multiplicity of $\mathscr{M}$}. The set of \emph{occurrences of $\mathscr{M}$} is defined by 
$$
O(\mathscr{M}):=\{(x,s): x\in M, 1\leq s\leq \mathsf{m}(x)\}\subset M\times \IN
$$

We first record the following generalization of the Isserlis-Wick theorem for Gaussian random vectors, where for a vector $(x_1,\dots,x_m)\in\IR^m$ and a multi-index $\alpha\in\IN^m_0$ we use the usual notation
$$
x^\alpha:=x_1^{\alpha_1}\cdots x_m^{\alpha_m}
$$
and for a natural number $d$ we set
$$
[d]:=\{1,\dots,d\}.
$$

Given $\alpha^{(1)}, \dots, \alpha^{(n)}\in\IN^m_0$ we define a multiset $\mathscr{M}(\alpha^{(1)}, \dots, \alpha^{(n)})$ by
$$
\mathscr{M}(\alpha^{(1)}, \dots, \alpha^{(n)}):=\left([n]\times [m], (k,i)\mapsto \alpha^{(k)}_i\right).
$$
Note that the set of occurrences of $\mathscr{M}(\alpha^{(1)}, \dots, \alpha^{(n)})$ is given by
\begin{align*}
O(\alpha^{(1)}, \dots, \alpha^{(n)})&:=O\big(\IMM(\alpha^{(1)}, \dots, \alpha^{(n)})\big)\\
&=\left\{(k,i,r): k\in [n], i\in [m], 1\leq r\leq\alpha^{(k)}_i\right\}\subset [n]\times [m]\times \IN,
\end{align*}
and that the cardinality of $O(\alpha^{(1)}, \dots, \alpha^{(n)})$ is given by 
$$
|O(\alpha^{(1)}, \dots, \alpha^{(n)})|=\sum^n_{k=1}\sum^m_{i=1}\alpha^{(k)}_i.
$$

With this notation at hand, we can now state:

\begin{proposition}\label{wick} For each $k=1,\dots,n$ let $Y^{(k)}=(Y^{(k)}_1,\dots, Y^{(k)}_m)\in\mathbb{R}^m$ be a zero mean Gaussian random vector, and let 
$$
Q(x)=\sum_{\alpha\in\IN^m_0}q_{\alpha} \, x^\alpha
$$
be a polynomial of $m$ variables with real coefficients. Then one has
\begin{align*}
&\mathbb{E}\left[\prod_{k=1}^nQ(Y^{(k)})\right]\\
&=\sum_{\substack{\alpha^{(1)}, \dots, \alpha^{(n)}\in\IN^m_0 \\|O(\alpha^{(1)}, \dots, \alpha^{(n)})|/2\,\in\IN}}\left(\prod_{k=1}^n q_{\alpha^{(k)}}\right)\,
\sum_{\substack{\text{\emph{pairings $\mathcal{P}$ of }}\\O(\alpha^{(1)}, \dots, \alpha^{(n)})}}\;
\prod_{\left\{\,(k,i,r)\, , \,(l,j,s)\,\right\}\in \mathcal{P}}\mathrm{Cov}\left(Y^{(k)}_{i},Y^{(l)}_{j}\right).
\end{align*}
\end{proposition}

\begin{proof} Expanding the product in the expectation gives
$$
\prod_{k=1}^nQ(Y^{(k)})=\sum_{\alpha^{(1)}, \dots, \alpha^{(n)}\in\IN^m_0 }\left(\prod_{k=1}^n q_{\alpha^{(k)}}\right)\,\left(\prod_{k=1}^n (Y^{(k)})^{\alpha^{(k)}}\right),
$$
so the result follows from the standard Isserlis' Theorem, as stated in the introduction.
\end{proof}

Given a matrix $A\in \mathbb{N}^{n\times l}_0$, we denote with\vspace{2mm}
\begin{itemize}
\setlength{\itemsep}{0.8em} 
\item $S(A)\in \mathbb{N}_0$ the sum of all entries of $A$,
\item $SC(A,i)\in \mathbb{N}_0$ the sum of all the entries of the $i$-th column of $A$, 
\item $SR(A,i)\in \mathbb{N}_0$ the sum of all the entries of the $i$-th row of $A$,
\item $A!\in  \mathbb{N}^{n\times l}_0$ the matrix received by applying the factorial function on each entry of $A$, 
\item $P(A)\in  \mathbb{N}_0$ the product of all entries, 
\item  $\tr(A)\in  \mathbb{N}_0$ the trace of $A$ in the case $l=n$.
\end{itemize}
\vspace{2mm}

Given another matrix $M\in \mathbb{N}^{n\times n}_0$, let
\begin{align*}
C(M,A) := \frac{P(A!)}{2^{\tr(M)}} 
\sum_{\substack{M_1,...,M_l\in\nat_0^{n\times n}:\\ 
M_1+...+M_l=M\\
\forall j\in[n]\forall k\in [l]:\\
SC(M_k,j)+SR(M_k,j)= A_{j,k}}} 
\prod_{i=1}^l \frac{1}{P(M_i!)}.
\end{align*}

Below, $C(M,A)$ will appear as a combinatorial factor for upper triangular $M$'s and vector $A$'s, where we count certain multigraphs. Regardless of this, we prove:

\begin{lemma} 
$C(M,A)$ is an integer.
\end{lemma}

\begin{proof} 
Fix $M_1,...,M_l\in\nat_0^{n\times n}$ such that the constraints
\[
SC(M_k,j)+SR(M_k,j) = A_{j,k} \quad \text{for all } j,k
\]
and
\[
M_1+\cdots+M_l = M
\]
hold. Write
\[
P(A!) = \prod_{j=1}^n \prod_{k=1}^l A_{j,k}!, 
\qquad 
P(M_i!) = \prod_{p,q=1}^n (M_i)_{p,q}!, \quad 
2^{\tr(M)} = 2^{\sum_{j=1}^n \sum_{k=1}^l (M_k)_{j,j}}.
\]
Thus, we have to show that
\[
z := \frac{\prod_{j=1}^n \prod_{k=1}^l A_{j,k}!}{2^{\sum_{j=1}^n \sum_{k=1}^l (M_k)_{j,j}}} \prod_{i=1}^l \frac{1}{\prod_{p,q=1}^n (M_i)_{p,q}!}
\]
is an integer. For this, let
\begin{equation*}
\begin{aligned}
E(c_1,...,c_n) &:= \frac{(c_1+...+c_n)!}{c_1!\cdots c_n!} \\
&= \binom{c_1+...+c_n}{c_1} \binom{c_2+...+c_n}{c_2} \cdots \binom{c_{n-1}+c_n}{c_{n-1}} \in \nat.
\end{aligned}
\end{equation*}
Then we have
\begin{equation*}
\begin{aligned}
z &= \frac{\prod_{j=1}^n \prod_{k=1}^l A_{j,k}!}{2^{\sum_{j=1}^n \sum_{k=1}^l (M_k)_{j,j}}} \prod_{i=1}^l \frac{1}{\prod_{p,q=1}^n (M_i)_{p,q}!} \\
&= \frac{1}{2^{\sum_{j=1}^n \sum_{k=1}^l (M_k)_{j,j}}} \prod_{i=1}^l \frac{\prod_{j=1}^n A_{j,i}!}{\prod_{p,q=1}^n (M_i)_{p,q}!} \\
&= \frac{1}{2^{\sum_{j=1}^n \sum_{k=1}^l (M_k)_{j,j}}} \prod_{i=1}^l \prod_{j=1}^n \frac{A_{j,i}!}{\prod_{q=1}^n (M_i)_{j,q}!} \\
&= \frac{1}{2^{\sum_{j=1}^n \sum_{k=1}^l (M_k)_{j,j}}} 
\prod_{i=1}^l \prod_{j=1}^n 
E((M_i)_{j,1},..., (M_i)_{j,n}, (M_i)_{1,j},..., (M_i)_{n,j}) \prod_{q=1}^n (M_i)_{q,j}! \\
&= \prod_{i=1}^l \prod_{j=1}^n 
\frac{1}{2^{(M_i)_{j,j}}} E((M_i)_{j,1},..., (M_i)_{j,n}, (M_i)_{1,j},..., (M_i)_{n,j}) \prod_{q=1}^n (M_i)_{q,j}!.
\end{aligned}
\end{equation*}

The function $E$ is symmetric and takes $(M_i)_{j,j}$ twice as a parameter. By reordering the parameters such that these two parameters are the last ones, we see that 
\[
E((M_i)_{j,1},..., (M_i)_{j,n}, (M_i)_{1,j},..., (M_i)_{n,j})
\] 
is divisible by 
\(\binom{(M_i)_{j,j}+(M_i)_{j,j}}{(M_i)_{j,j}}\). Considering the other parts of the last term, it suffices to show that
\[
\binom{(M_i)_{j,j}+(M_i)_{j,j}}{(M_i)_{j,j}} \frac{(M_i)_{j,j}!}{2^{(M_i)_{j,j}}}
\]
is an integer. Using the notation \(y=(M_i)_{j,j}\), we have
\begin{equation*}
\begin{aligned}
\binom{y+y}{y} \frac{y!}{2^y} &= \frac{(2y)!}{y! 2^y} \\
&= \frac{(1\cdot 3 \cdots (2y-1)) \cdot (2\cdot 4 \cdots 2y)}{y! 2^y} \\
&= \frac{(1\cdot 3 \cdots (2y-1)) \cdot (1\cdot 2 \cdots y)}{y!} \\
&= 1\cdot 3 \cdots (2y-1) \in \nat,
\end{aligned}
\end{equation*}
which proves that $z \in \mathbb{N}$.
\end{proof}

An \emph{undirected finite multigraph with loops} (short: \emph{multigraph}) $\Gamma$ is a pair $\Gamma=(V,E)$, where $V$ is any countable set and $E=(X,h)$ is a multiset with universe
    $$
        X=\{\{v,w\}: v,w\in V\}.
    $$
    Above, $V$ is called the \emph{vertices} and $h$ the \emph{edge-muliplicity function} of $\Gamma$.   
\vspace{1mm}

Note that $X$ also contains all singletons $\{v\}$ for $v\in V$. For $v\in V$ we also define the \emph{degree} as
    $$
        \deg(v)=2h(\{v\})+\sum_{w\in V\backslash\{v\}}h(\{v,w\}),
    $$
so loops count twice to the degree of a vertex. Vertices $v$, $w$ are interpreted to have $h(\{v,w\})$ many edges joining them.
\vspace{1mm}

The multigraph $\Gamma$ is called \emph{connected}, if for all vertices $a\neq b$ there are vertices $v_1,\dots,v_n$ such that 
$$
h(\{a,v_1\})\geq 1,\quad h(\{v_i,v_{i+1}\})\geq 1\quad\text{for all $i=1,\dots,n-1$},\quad h(\{v_n,b\})\geq 1.
$$

A \emph{sub-multigraph} of $\Gamma$ is a multigraph $\Gamma'=(V',E')$ with $V'\subset V$ whose edge-multiplicity function is the restriction of $h$ to the universe of $E'$.\vspace{1mm}

A \emph{connected component of $\Gamma$} is a connected sub-multigraph of $\Gamma$ which is not a sub-multigraph of any other connected sub-multigraph of $\Gamma$.  We denote 
$$
\mathscr{C}(\Gamma):=\{\text{connected components of $\Gamma$}\}.
$$

\vspace{1mm}

If we denote the vertices of $\Gamma$ with 
$$
V=\{v_{1},...,v_{l}\},
$$
then given a symmetric function
    $$
        f:[0,1]\times [0,1]\longrightarrow  \mathbb{R}, 
    $$
we define 
    $$
        \int_{\Gamma}  f:= \int_{[0,1]^l} \prod_{(\{v_a,v_b\},j)\in O(E)}f(s_a,s_b)\ \Id s_{1}\cdots\Id s_{l}.
    $$
This definition is independent of the order chosen on the vertices $V$, since we can just reorder the integration variables and it is independent of the order chosen on the elements of the universe of $E$, since $f$ is symmetric.\\
We define the natural number $C(\Gamma)$ by 
$$
C(\Gamma):=C(M^{\Gamma},A^{\Gamma}),
$$
with $M^{\Gamma}\in \mathbb{N}^{l\times l}_0$ the upper triangular adjacency matrix of $\Gamma$, so
    $$
    M^{\Gamma}_{ij}:=\begin{cases} h(\{v_i,v_j\}),\quad\text{if $i< j$}\\h(\{v_i\}),\quad\text{if $i=j$}\\
    0,\quad\text{if $i> j$}\end{cases} ,
    $$
   and
    $$
    A^{\Gamma}:=\pmat{SC(M^{\Gamma},1)+SR(M^{\Gamma},1)\\\vdots\\SC(M^{\Gamma},l)+SR(M^{\Gamma},l)}=\pmat{\deg(v_1)\\\vdots\\\deg(v_l)}\in \mathbb{N}^{l\times 1}_0.
    $$
From the definition of $C(\Gamma)$ we immediately have
\begin{align}\label{edya}
C(\Gamma)=\frac{P(A^\Gamma!)}{2^{\tr(M^\Gamma)}} \frac{1}{P(M^\Gamma!)}.
\end{align}
This definition is also independent of the order chosen on the vertices $V$, because a permutation of the vertices will just permute the elements of $A^\Gamma$ and permute the entries of $M^\Gamma$, where each diagonal entry of $M^\Gamma$ will stay on the diagonal, so all terms in the above equation stay the same.\\
When we have two multigraphs $\Gamma_1=(V,E_1)$ and $\Gamma_2=(V,E_2)$ with the same vertex set, we can define their \emph{sum} as follows: when $E_1=(X,h_1)$ and $E_2=(X,h_2)$, we define
    $$
        \Gamma_1+\Gamma_2:=(V,E_3)\qquad\t{with}\qquad E_3=(X,h_1+h_2),
    $$
    where $h_1+h_2$ is the componentwise sum of the maps $h_1$ and $h_2$. Clearly, this sum is commutative and associative, so we can omit brackets when dealing with sums of more than two graphs.\vspace{1mm}


\begin{theorem}\label{propwithgraphs}
    Assume that
    $$
        f:[0,1]\times [0,1]\longrightarrow  \mathbb{R}
    $$
    is a continuous symmetric function, and that $\alpha^{(1)}, \dots, \alpha^{(n)}\in\nat_0^m$
    are such that $\alpha_j^{(1)}+ \dots+ \alpha_j^{(n)}$ is even for all $j\in[m]$. Fix also a finite set $\{v_1,\dots, v_n\}$ of vertices. Then one has
    \begin{align*}
        &\int_{[0,1]^n}
        \sum_{\substack{\text{\emph{pairings} }\mathcal{P}\\\text{of }O(\alpha^{(1)},\dots,\alpha^{(n)})}}
        \prod_{\{(k,i,r),(l,j,s)\}\in\mathcal{P}} \delta_{ij}\, f(s_k,s_l)\, \Id s_1\cdots \Id s_n\\
        =&\sum_{\substack{\Gamma_1=(\{v_1,...,v_n\},E_1)\\\forall j\in[n]:\;\deg(v_j)=\alp_1^{(j)}}}\dots\sum_{\substack{\Gamma_m=(\{v_1,...,v_n\},E_m)\\\forall j\in[n]:\;\deg(v_j)=\alp_m^{(j)}}}\left(\prod_{q=1}^m C(\Gamma_q)\right)\left(\prod_{\Lambda\in\mathscr{C}(\Gamma_1+\dots\Gamma_m)}\int_{\Lambda} f\right).
    \end{align*}
\end{theorem}

\begin{proof}
    We define 
    $$
        B_q:=O(([n],k\mapsto \alp_q^{(k)}))\quad\text{for all $q\in[m]$.}
    $$
    On the left-hand side of the claimed equation, all pairings that contain an element of the form $\{(k,i,r),(l,j,t)\}$ with $i\neq j$ are irrelevant since the product will be 0, so we can split up the pairings:
    \begin{align*}
        &\int_{[0,1]^n}\sum_{\substack{\t{pairings $\cp$}\\\t{of }O(\alp^{(1)},...,\alp^{(n)})}}\prod_{\{(k,i,r),(l,j,t)\}\in\cp}\delta_{ij}\;f(s_k,s_l)\;\mathrm{d}s_1\cdots\mathrm{d}s_n\\
        =&\int_{[0,1]^n}\sum_{\substack{\t{pairings $\cp_1$}\\\t{of }B_1}}\cdots \sum_{\substack{\t{pairings $\cp_m$}\\\t{of }B_m}}\\&\qquad\qquad\left(\prod_{\{(k_1,r_1),(l_1,t_1)\}\in\cp_1}f(s_{k_1},s_{l_1})\right)\cdots\left(\prod_{\{(k_m,r_m),(l_m,t_m)\}\in\cp_m}f(s_{k_m},s_{l_m})\right)\;\mathrm{d}s_1\cdots\mathrm{d}s_n.
    \end{align*}
    For any $q\in[m]$, a pairing $\cp_q$ of $B_q$ induces a multigraph $\Gamma_q=(V,E_q)$ with
    $$
        V=\{v_1,...,v_n\}
    $$
    and
    $$
        \deg(v_j)=\alp^{(j)}\text{ for all $v_j\in V$,}
    $$
    by interpreting each pair $\{(k,r),(l,s)\}\in\cp$ as an edge between $v_k$ and $v_l$. Any such multigraph, however, is (usually) induced by multiple pairings, say by $\tilde{C}(\Gamma_q)$ many different pairings. By interpreting the pairings as multigraphs, we get
    \begin{align*}
        &\int_{[0,1]^n}\sum_{\substack{\t{pairings $\cp_1$}\\\t{of }B_1}}\cdots \sum_{\substack{\t{pairings $\cp_m$}\\\t{of }B_m}}\\&\qquad\qquad\left(\prod_{\{(k_1,r_1),(l_1,t_1)\}\in\cp_1}f(s_{k_1},s_{l_1})\right)\cdots\left(\prod_{\{(k_m,r_m),(l_m,t_m)\}\in\cp_m}f(s_{k_m},s_{l_m})\right)\;\mathrm{d}s_1\cdots\mathrm{d}s_n\\
        =&\sum_{\substack{\Gamma_1=(\{v_1,...,v_n\},E_1)\\\forall j\in[n]:\;\deg(v_j)=\alp_1^{(j)}}}\dots\sum_{\substack{\Gamma_m=(\{v_1,...,v_n\},E_m)\\\forall j\in[n]:\;\deg(v_j)=\alp_m^{(j)}}}\tilde{C}(\Gamma_1)\cdots \tilde{C}(\Gamma_m)\int_{[0,1]^n}\\&\qquad\qquad\left(\prod_{(\{v_a,v_b\},j)\in O(E_1)}f(s_a,s_b)\right)\cdots\left(\prod_{(\{v_a,v_b\},j)\in O(E_m)}f(s_a,s_b)\right)\, \Id s_1\cdots \Id s_n\\
        =&\sum_{\substack{\Gamma_1=(\{v_1,...,v_n\},E_1)\\\forall j\in[n]:\;\deg(v_j)=\alp_1^{(j)}}}\dots\sum_{\substack{\Gamma_m=(\{v_1,...,v_n\},E_m)\\\forall j\in[n]:\;\deg(v_j)=\alp_m^{(j)}}}\left(\prod_{q=1}^m \tilde{C}(\Gamma_q)\right)\int_{\Gamma_1+\dots+\Gamma_m}  f\\
        =&\sum_{\substack{\Gamma_1=(\{v_1,...,v_n\},E_1)\\\forall j\in[n]:\;\deg(v_j)=\alp_1^{(j)}}}\dots\sum_{\substack{\Gamma_m=(\{v_1,...,v_n\},E_m)\\\forall j\in[n]:\;\deg(v_j)=\alp_m^{(j)}}}\left(\prod_{q=1}^m \tilde{C}(\Gamma_q)\right)\prod_{\Lambda\in \mathscr{C}(\Gamma_1+\dots+\Gamma_m)}\int_{\Lambda}  f.
    \end{align*}
    It remains to show the asserted formula formula for $\tilde{C}(\Gamma_q)$:
    \begin{align*}
        \tilde{C}(\Gamma_q)=&|\t{pairings }\cp\t{ of }B_q\t{ s.t. }\forall l\in[n]\forall k\in[l]:\\
        &\qquad\t{There are exactly }h(\{v_k,v_l\})\t{ many pairings of the form }\{(k,r),(l,s)\}\in\cp|\\
        =&|\t{pairings }\cp\t{ of }B_q\t{ s.t. }\forall l\in[n]\forall k\in[l]:\\
        &\qquad\t{There are exactly }M_{k,l}^{\Gamma_q}\t{ many pairings of the form }\{(k,r),(l,s)\}\in\cp|\\
        =&:\tilde C(M^{\Gamma_q},A^{\Gamma_q})
    \end{align*}
    Now we only need to show that $\tilde C(M^{\Gamma_q},A^{\Gamma_q})=C(M^{\Gamma_q},A^{\Gamma_q})$.
    The term
    $$
        \prod_{\{(k_q,r_q),(l_q,t_q)\}\in\cp_q}f(s_{k_q},s_{l_q})
    $$
    contains ''$s_i$'' exactly $\alpha_q^{(i)}$ times, for any pairing $\cp_q$ of $B_q$.\\
    The term
    $$
        \prod_{(\{v_a,v_b\},j)\in O(E_q)}f(s_a,s_b)
    $$
    contains ''$s_i$'' exactly $\deg(v_i)=SC(M^{\Gamma_q},i)+SR(M^{\Gamma_q},i)$ times.\\
    We can clearly find a pairing $\cp_a$ of $B_q$ with
    $$
        \prod_{\{(k,r),(l,t)\}\in\cp_a}\;f(s_k,s_l)=\prod_{1\leq k\leq l\leq n} f(s_k,s_l)^{M_{k,l}^{\Gamma_q}}.
    $$
    All other pairings with this property (where we read the equation above ''symbolically''; it could hold for other pairings if $f$ is chosen ''luckily'') can be defined in the following way:
    $$
        \cp_{\sigma_1,...,\sigma_n}=\{\{(k,\sigma_k(r)),(l,\sigma_l(t))\}\t{ for } \{(k,r),(l,t)\}\in\cp_a\},
    $$
    where $\sigma_j$ is a permutation on $\alp_q^{(j)}$ elements for all $j\in[n]$. So $\tilde C(M^{\Gamma_q},A^{\Gamma_q})$ should be equal to $\alp_q^{(1)}!\c...\c\alp_q^{(n)}!$, but there are 2 subtleties which we need to address.\\
    When 
    $$
    \{(k,r_1),(l,t_1)\},\{(k,r_2),(l,t_2)\},...,\{(k,r_q),(l,t_q)\}\in\cp_a
    $$
    for pairwise distinct $t_i$ and pairwise distinct $r_i$, then there are $q!$ ways to choose $\sigma_k$ and $\sigma_l$ that result in the same $\cp_{\sigma_1,...,\sigma_n}$, so we need to divide by $q!$. Since the $q$'s are exactly the entries of $M^{\Gamma_q}$, we need to divide by $P(M^{\Gamma_q}!)$.\\
    When 
    $$
    \{(k,r_1),(k,t_1)\}\in\cp_a
    $$
    for $r_1\neq t_1$, we can replace $\sigma_k$ by $\sigma_k\circ\tau$, where $\tau$ is the transposition of $r_1$ and $t_1$ to get the same $\cp_{\sigma_1,...,\sigma_n}$. So we need to divide by 2. This happens exactly $\tr(M^{\Gamma_q})$ times, so we need to divide by $2^{\tr(M^{\Gamma_q})}$. So we get
    \begin{align*}
        \tilde{C}(M^{\Gamma_q},A^{\Gamma_q})=&\>\f{\alp_q^{(1)}!\c...\c\alp_q^{(n)}!}{P(M^{\Gamma_q}!)2^{\tr(M^{\Gamma_q})}}\\
        =&\>C(M^{\Gamma_q},A^{\Gamma_q}),
    \end{align*}
    completing the proof.
\end{proof}

We immediately obtain:

\begin{theorem}\label{prop1gen} Assume 
$$
X(s)=(X_1(s),\dots,X_m(s))\in \mathbb{R}^m,\quad s\in [0,1], 
$$
is a zero mean Gaussian process with a jointly continuous covariance matrix of the form 
$$
\mathrm{Cov}(X_i(s),X_j(t))=f(s,t)\delta_{ij},\quad\text{for some continuous function $f:[0,1]\times [0,1]\to\mathbb{R}$,}
$$
and assume
$$
Q(x)=\sum_{\alpha\in\IN^m_0}q_{\alpha} \, x^\alpha
$$
is a polynomial of $m$ variables with real coefficients. Fix a finite set $\{v_1,...,v_n\}$ of vertices. Then one has
\begin{align*}
&\int_{\{0\leq s_1\leq...\leq s_n\leq 1\}}\mathbb{E}\left[\prod_{k=1}^nQ(X(s_k))\right]\;\mathrm{d}s_1\cdots\mathrm{d}s_n\\
&=\f1{n!}\sum_{\substack{\alpha^{(1)}, \dots, \alpha^{(n)}\in\IN^m_0 \\|O(\alpha^{(1)}, \dots, \alpha^{(n)})|/2\,\in\IN}}\left(\prod_{k=1}^n q_{\alpha^{(k)}}\right)\times\\
&\quad \times \sum_{\substack{\Gamma_1=(\{v_1,...,v_n\},E_1)\\\forall j\in[n]:\;\deg(v_j)=\alp_1^{(j)}}}\dots\sum_{\substack{\Gamma_m=(\{v_1,...,v_n\},E_m)\\\forall j\in[n]:\;\deg(v_j)=\alp_m^{(j)}}}\left(\prod_{q=1}^m C(\Gamma_q)\right)\left(\prod_{\Lambda\in \mathscr{C}(\Gamma_1+\cdots+\Gamma_m)}\int_{\Lambda}  f\right).
\end{align*}
\end{theorem}

\begin{proof} This follows from combining
\begin{align*}
&\int_{\{0\leq s_1\leq...\leq s_n\leq 1\}}\mathbb{E}\left[\prod_{k=1}^nQ(X(s_k))\right]\;\mathrm{d}s_1\cdots\mathrm{d}s_n=\frac{1}{n!}\int_{[0,1]^n}\mathbb{E}\left[\prod_{k=1}^nQ(X(s_k))\right]\;\mathrm{d}s_1\cdots\mathrm{d}s_n
\end{align*}
with Theorem \ref{propwithgraphs}  and Proposition \ref{wick}.

\end{proof}

\begin{example}\label{edxaa} We take $c\in{\real^{m\times 1}}$, $D\in{\real^{m\times m}}$ symmetric and the polynomial
\begin{equation*}
    Q(x)=(x,Dx)+(c,x)=:\sum_{\alpha\in\mathbb{N}^m_0}q_\alpha x^\alpha.
\end{equation*}
Then we have 
$$
q_{e_i}=c_i, \quad q_{2e_i}=D_{i,i},\quad q_{e_i+e_j}=2D_{i,j}\quad\text{if $i\neq j$},
$$
with all other $q_\alpha=0$. \\
$\mathbf{n=1}$: We fix a set with one element $\{v\}$ and consider the expression
\begin{align*}
\Psi_1[f]:=\sum_{\substack{\alpha\in\IN^m_0 \\|O(\alpha)|/2\,\in\IN}} q_{\alpha}\sum_{\substack{\Gamma_1=(\{v\},E_1)\\\deg(v)=\alp_1}}\dots\sum_{\substack{\Gamma_m=(\{v\},E_m)\\\deg(v)=\alp_m}}\left(\prod_{q=1}^m C(\Gamma_q)\right)\left(\prod_{\Lambda\in \mathscr{C}(\Gamma_1+\dots\Gamma_m)}\int_{\Lambda}  f\right).
\end{align*}
The only $q_\alpha$ with $|O(\alpha)|$ even and a nontrivial $\Gamma_l$-sum is given by $q_{2e_i}=D_{i,i}$. For this choice of $\alpha$, the only nontrivial sum in $\Gamma_l$-sum is the $i$-th one, which contains precisely the graph
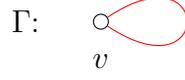
\begin{center}
\begin{tikzpicture}
    \node[circle, draw, inner sep=2pt] (Ab2) at (0, 2) {};
   \node[below=5pt of Ab2] {$v$};
    \draw[-, red] (Ab2) .. controls (1.5,3) and (1.5, 1) .. (Ab2);
    \node at (-1,2) {$\Gamma$:};
\end{tikzpicture}
\end{center}
and we get the contribution
$$
D_{i,i}C(\Gamma)\int_{\Gamma}f=D_{i,i}\int_{[0,1]} f(s,s) \Id s.
$$
and so
$$
\Psi_1[f]=\sum_i D_{i,i}\int_{[0,1]} f(s,s) \Id s.
$$

$\mathbf{n=2}$: We fix a set with two elements $\{v,w\}$ and consider the expression
\begin{align*}
\Psi_2[f]:=&\f1{2}\sum_{\substack{\alpha,  \beta\in\IN^m_0 \\|O(\alpha, \beta)|/2\,\in\IN}} q_{\alpha}q_{\beta}\sum_{\substack{\Gamma_1=(\{v,w\},E_1)\\\deg(v)=\alp_1,\deg(w)=\beta_1}}\times\dots\\&\times\sum_{\substack{\Gamma_m=(\{v,w\},E_m)\\\deg(v)=\alp_m,\deg(w)=\beta_m}}\left(\prod_{q=1}^m C(\Gamma_q)\right)\left(\prod_{\Lambda\in \mathscr{C}(\Gamma_1+\dots\Gamma_m)}\int_{\Lambda}  f\right).
\end{align*}
The $q_\alpha q_ \beta$ with $|O(\alpha, \beta)|$ even and with some nontrivial $\Gamma_l$-sum are given by
\begin{align*}
&q_{e_i}q_{e_i}=c_i^2,\quad q_{2e_i}q_{2e_i}=D_{i,i}^2,\quad q_{e_i+e_j}q_{e_i+e_j}=4D_{i,j}^2\quad\text{with $i\neq j$}, \\
&q_{2e_i}q_{2e_j}=D_{i,i}D_{j,j}\quad\text{with $i\neq j$}.
\end{align*}
For the first choice of $\alpha,\beta$, the only nontrivial sum in $\sum_{\Gamma_1}\cdots\sum_{\Gamma_m}$ is the $i$-th one, which has precisely one summand given by the graph 
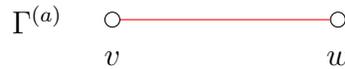
\begin{center}
\begin{tikzpicture}
    \node[circle, draw, inner sep=2pt] (Aa) at (0, 6) {};
    \node[circle, draw, inner sep=2pt] (Ba) at (3, 6) {};
    \node[below=5pt of Aa] {$v$};
    \node[below=5pt of Ba] {$w$};
    \draw[-, red] (Aa) -- (Ba);
    \node at (-1,6) {$\Gamma^{(a)}$};
\end{tikzpicture}
\end{center}
Here we get the contribution
$$
\frac{1}{2}c_i^2C(\Gamma^{(a)})\int_{\Gamma^{(a)}}f=\frac{1}{2}c_i^2\int_{[0,1]^2} f(s_1,s_2)\Id s_1 \Id s_2.
$$
For the second choice of $\alpha,\beta$, again all sums are trivial except the $i$-th one, which now contains the two graphs

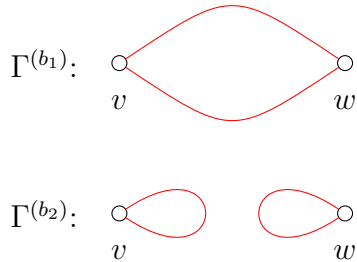
\begin{center}
\begin{tikzpicture}
    \node[circle, draw, inner sep=2pt] (Ab) at (0, 4) {};
    \node[circle, draw, inner sep=2pt] (Bb) at (3, 4) {};
    \node[below=5pt of Ab] {$v$};
    \node[below=5pt of Bb] {$w$};
    \node[right=5pt of Bb] {};
    \draw[-, red] (Ab) .. controls (1.5, 5) and (1.5, 5) .. (Bb);
    \draw[-, red] (Ab) .. controls (1.5, 3) and (1.5, 3) .. (Bb);
    \node at (-1,4) {$\Gamma^{(b_1)}$:};
    \node[circle, draw, inner sep=2pt] (Ab2) at (0, 2) {};
    \node[circle, draw, inner sep=2pt] (Bb2) at (3, 2) {};
    \node[below=5pt of Ab2] {$v$};
    \node[below=5pt of Bb2] {$w$};
    \node[right=5pt of Bb2] {};
    \draw[-, red] (Ab2) .. controls (1.5,3) and (1.5, 1) .. (Ab2);
    \draw[-, red] (Bb2) .. controls (1.5,3) and (1.5, 1) .. (Bb2);
    \node at (-1,2) {$\Gamma^{(b_2)}$:};
\end{tikzpicture}
\end{center}
We get the contributions 
$$
\frac{1}{2}D_{i,i}^2C(\Gamma^{(b_1)})\int_{\Gamma^{(b_1)}}f= D_{i,i}^2\int_{[0,1]^2} f(s_1,s_2)^2 \Id s_1\Id s_2,
$$
and
$$
\frac{1}{2}D_{i,i}^2C(\Gamma^{(b_2)})\int_{\Gamma^{(b_2)}}f= \frac{1}{2}D_{i,i}^2\left(\int_{[0,1]^2} f(s,s) \Id s\right)^2.
$$
The next choice of $\alpha,\beta$, which corresponds to $q_{e_i+e_j}q_{e_i+e_j}$, leads to two nontrivial sums, the $i$-th one and the $j$-th one. Each of these sums contains precisely the graph $\Gamma^{(a)}$, and their sum $2\Gamma^{(a)}$ is the connected graph $\Gamma^{(b_1)}$, so that the overall contribution is
$$
\frac{1}{2}4D_{i,j}^2C(\Gamma^{(a)})^2\int_{\Gamma^{(b_1)}}f=2D_{i,j}^2\int_{[0,1]^2} f(s_1,s_2)^2 \Id s_1\Id s_2.
$$
The last choice of $\alpha,\beta$, which corresponds to $q_{2e_i}q_{2e_j}$, leads again to two nontrivial sums, the $i$-th one and the $j$-th one. Each of these sums contains precisely one summand, the first one being 
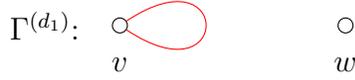
\begin{center}
\begin{tikzpicture}
\node[circle, draw, inner sep=2pt] (Ad) at (0, -2) {};
    \node[circle, draw, inner sep=2pt] (Bd) at (3, -2) {};
    \node[below=5pt of Ad] {$v$};
    \node[below=5pt of Bd] {$w$};
    \node[right=5pt of Bd] {};
    \draw[-, red] (Ad) .. controls (1.5,-1) and (1.5, -3) .. (Ad);
    \node at (-1,-2) {$\Gamma^{(d_1)}$:};
\end{tikzpicture}
\end{center}
the second one being
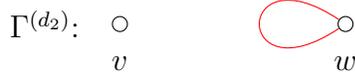
\begin{center}
\begin{tikzpicture}
  \node[circle, draw, inner sep=2pt] (Ad) at (0, -2) {};
    \node[circle, draw, inner sep=2pt] (Bd) at (3, -2) {};
    \node[below=5pt of Ad] {$v$};
    \node[below=5pt of Bd] {$w$};
    \node[right=5pt of Bd] {};
    \draw[red] (Bd) .. controls (1.5,-1) and (1.5, -3) .. (Bd);
    \node at (-1,-2) {$\Gamma^{(d_2)}$:};
\end{tikzpicture}
\end{center}
and their sum $\Gamma^{(d_1)}+\Gamma^{(d_2)}$ equals $\Gamma^{(b_2)}$. Thus, the overall contribution is
$$
\frac{1}{2}D_{i,i}D_{j,j}C(\Gamma^{(d_1)})C(\Gamma^{(d_2)})\int_{\Gamma^{(b_2)}}f=\frac{1}{2}D_{i,i}D_{j,j}\left(\int_{[0,1]^2} f(s,s) \Id s\right)^2.
$$
Putting everything together we arrive at
\begin{align*}
\Psi_2[f]&=\frac{1}{2} \sum_{i}c_i^2 \int_{[0,1]^2} f(s_1,s_2)\Id s_1 \Id s_2+ \sum_i D_{i,i}^2 \int_{[0,1]^2} f(s_1,s_2)^2 \Id s_1\Id s_2\\
&+\frac{1}{2}  \sum_i D_{i,i}^2 \left(\int_{[0,1]} f(s,s)\Id s\right)^2 +2 \sum_{i\neq j} D_{i,j}^2 \int_{[0,1]^2} f(s_1,s_2)^2 \Id s_1\Id s_2\\
&+ \frac{1}{2}\sum_{i\neq j} D_{i,i}D_{j,j} \left(\int_{[0,1]} f(s,s)\Id s\right)^2.
\end{align*}
$\mathbf{n=3}$: We fix a set with three elements $\{v,w,x\}$ and consider the expression
\begin{align*}
\Psi_3[f]:=&\f1{6}\sum_{\substack{\alpha,  \beta,\gamma\in\IN^m_0 \\|O(\alpha, \beta,\gamma)|/2\,\in\IN}} q_{\alpha}q_{\beta}q_{\gamma}\sum_{\substack{\Gamma_1=(\{v,w,x\},E_1)\\\deg(v)=\alp_1,\deg(w)=\beta_1,\deg(x)=\gamma_1}}\dots\times\\
&\times\dots\sum_{\substack{\Gamma_m=(\{v,w,x\},E_m)\\\deg(v)=\alp_m,\deg(w)=\beta_m,\deg(x)=\gamma_m}}\left(\prod_{q=1}^m C(\Gamma_q)\right)\left(\prod_{\Lambda\in \mathscr{C}(\Gamma_1+\dots\Gamma_m)}\int_{\Lambda}  f\right)
\end{align*}
The $q_\alpha q_\beta q_\gamma$ with $|O(\alpha, \beta,\gamma)|$ even with some nontrivial $\Gamma_l$-sum are given by (always with $i\neq j\neq k\neq i$)
\begin{align*}
    &q_{e_i}q_{2e_i}q_{e_i}=c_i^2D_{i,i},\quad q_{e_j}q_{2e_i}q_{e_j}=D_{i,i}c_j^2,\quad q_{e_i}q_{e_i+e_j}q_{e_j}=2c_ic_jD_{i,j}, \\
    &q_{e_i+e_j}q_{2e_i}q_{e_i+e_j}=4D_{i,j}^2D_{i,i},\quad q_{e_i+e_j}q_{2e_k}q_{e_i+e_j}=4D_{i,j}^2D_{k,k},\quad q_{2e_i}q_{2e_j}q_{2e_i}=D_{i,i}^2D_{j,j}\\
    &q_{2e_i}q_{2e_i}q_{2e_i}=D_{i,i}^3,\quad q_{2e_i}q_{2e_j}q_{2e_k}=D_{i,i}D_{j,j}D_{k,k},\quad q_{e_i+e_j}q_{e_i+e_k}q_{e_j+e_k}=8D_{i,j}D_{i,k}D_{j,k},
\end{align*}
where all the terms in the first two rows also have variants with exchanged $\alpha, \beta, \gamma$ that cannot be reached just by changing $i,j,k$. The corresponding terms need to be multiplied by 3.\\
The choice of $\alpha,\beta,\gamma$ which corresponds to $q_{e_i}q_{2e_i}q_{e_i}$ leads to only trivial sums except the $i$-th one, where we get two graphs
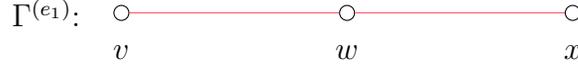
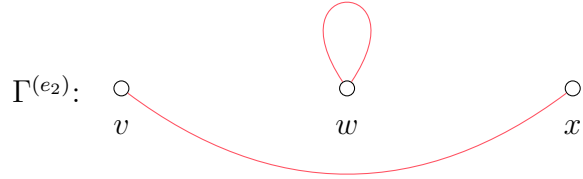
\threevertexgraph000100010{\Gamma^{(e_1)}}
and
\threevertexgraph010001000{\Gamma^{(e_2)}}
The corresponding contributions are
$$
    \frac{3}{6}c_i^2D_{i,i}C(\Gamma^{(e_1)})\int_{\Gamma^{(e_1)}}f=c_i^2D_{i,i}\int_{[0,1]^3} f(s_1,s_2)f(s_2,s_3) \Id s_1\Id s_2\Id s_3
$$
and
$$
    \frac{3}{6}c_i^2D_{i,i} C(\Gamma^{(e_2)})\int_{\Gamma^{(e_2)}}f=\frac{1}{2}c_i^2D_{i,i}\int_{[0,1]^2} f(s_1,s_2) \Id s_1\Id s_2\int_{[0,1]} f(s,s) \Id s.
$$
The choice of $\alpha,\beta,\gamma$ which corresponds to $q_{e_j}q_{2e_i}q_{e_j}$ leads to only trivial sums except the $i$-th and $j$-th one, each of which only contains one possible graph. We draw the graphs in the same diagram, where red edges are in $\Gamma_i$ and blue edges are in $\Gamma_j$:
\threevertexgraph010002000{\Gamma^{(f)}}
The corresponding contribution is
$$
    \f36 D_{i,i}c_j^2C(\Gamma^{(f)}_{\t{red}})C(\Gamma^{(f)}_{\t{blue}})\int_{\Gamma^{(f)}}f=\f12 D_{i,i}c_j^2\int_{[0,1]} f(s,s) \Id s\int_{[0,1]^2} f(s_1,s_2) \Id s_1\Id s_2.
$$
The choice of $\alpha,\beta,\gamma$ which corresponds to $q_{e_i}q_{e_i+e_j}q_{e_j}$ leads to only trivial sums except the $i$-th and $j$-th one, each of which only contains one possible graph:
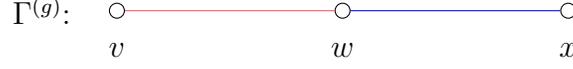
\threevertexgraph000100020{\Gamma^{(g)}}
The corresponding contribution is
$$
    \f36 2c_ic_jD_{i,j}C(\Gamma^{(g)}_{\t{red}})C(\Gamma^{(g)}_{\t{blue}})\int_{\Gamma^{(g)}}f=c_ic_jD_{i,j}\int_{[0,1]^3} f(s_1,s_2)f(s_2,s_3) \Id s_1\Id s_2\Id s_3.
$$
The choice of $\alpha,\beta,\gamma$ which corresponds to $q_{e_i+e_j}q_{2e_i}q_{e_i+e_j}$ leads to only trivial sums except the $i$-th and $j$-th one, where the possible graphs are
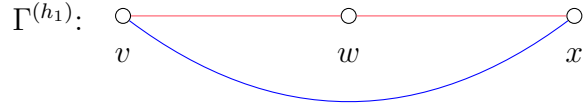
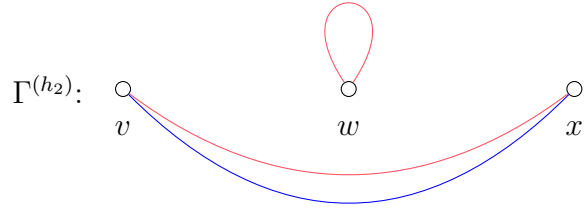
\threevertexgraph000102010{\Gamma^{(h_1)}}
and
\threevertexgraph010001200{\Gamma^{(h_2)}}
The corresponding contributions are
$$
    \f36 4D_{i,j}^2D_{i,i}C(\Gamma^{(h_1)}_{\t{red}})C(\Gamma^{(h_1)}_{\t{blue}})\int_{\Gamma^{(h_1)}}f=4D_{i,j}^2D_{i,i}\int_{[0,1]^3} f(s_1,s_2)f(s_2,s_3)f(s_3,s_1) \Id s_1\Id s_2\Id s_3.
$$
and
$$
    \f36 4D_{i,j}^2D_{i,i}C(\Gamma^{(h_2)}_{\t{red}})C(\Gamma^{(h_2)}_{\t{blue}})\int_{\Gamma^{(h_2)}}f=2D_{i,j}^2D_{i,i}\int_{[0,1]^2} f(s_1,s_2)^2 \Id s_1\Id s_2\int_{[0,1]} f(s,s) \Id s.
$$
\end{example}
The choice of $\alpha,\beta,\gamma$ which corresponds to $q_{e_i+e_j}q_{2e_k}q_{e_i+e_j}$ leads to only trivial sums except the $i$-th, $j$-th and $k$-th one. The edges of the graph $\Gamma_k$ get a green color:
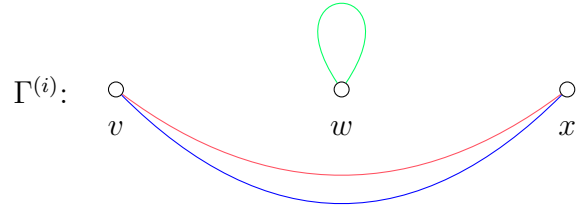
\threevertexgraph030001200{\Gamma^{(i)}}
The corresponding contribution is
$$
    \f36 4D_{i,j}^2D_{k,k}C(\Gamma^{(i)}_{\t{red}})C(\Gamma^{(i)}_{\t{blue}})C(\Gamma^{(i)}_{\t{green}})\int_{\Gamma^{(i)}}f=2D_{i,j}^2D_{k,k}\int_{[0,1]^2} f(s_1,s_2)^2 \Id s_1\Id s_2\int_{[0,1]} f(s,s) \Id s.
$$
The choice of $\alpha,\beta,\gamma$ which corresponds to $q_{2e_i}q_{2e_j}q_{2e_i}$ leads to only trivial sums except the $i$-th and $j$-th one, where the possible graphs are
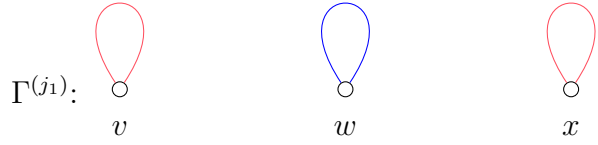
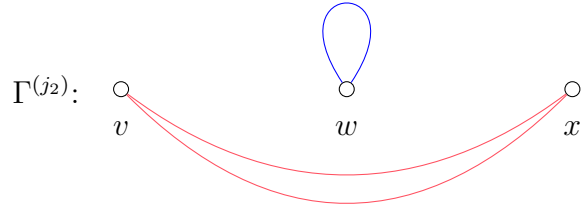
\threevertexgraph121000000{\Gamma^{(j_1)}}
and
\threevertexgraph020001100{\Gamma^{(j_2)}}
The corresponding contributions are
$$
    \f36 D_{i,i}^2D_{j,j}C(\Gamma^{(j_1)}_{\t{red}})C(\Gamma^{(j_1)}_{\t{blue}})C(\Gamma^{(j_1)}_{\t{green}})\int_{\Gamma^{(j_1)}}f=\f12D_{i,i}^2D_{j,j}\lb\int_{[0,1]}f(s,s)\Id s\rb^3
$$
and
$$
    \f36 D_{i,i}^2D_{j,j}C(\Gamma^{(j_2)}_{\t{red}})C(\Gamma^{(j_2)}_{\t{blue}})C(\Gamma^{(j_2)}_{\t{green}})\int_{\Gamma^{(j_2)}}f=D_{i,i}^2D_{j,j}\int_{[0,1]^2} f(s_1,s_2)^2 \Id s_1\Id s_2\int_{[0,1]} f(s,s) \Id s.
$$
The choice of $\alpha,\beta,\gamma$ which corresponds to $q_{2e_i}q_{2e_i}q_{2e_i}$ leads to only trivial sums except the $i$-th one, where we have five possible graphs:
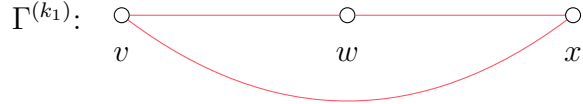
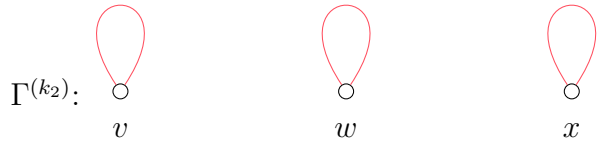
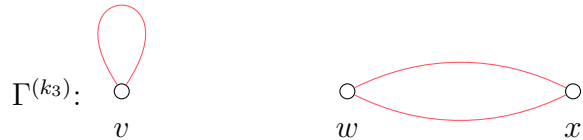
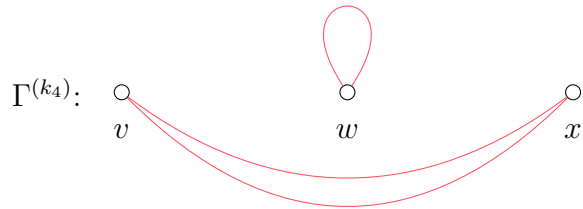
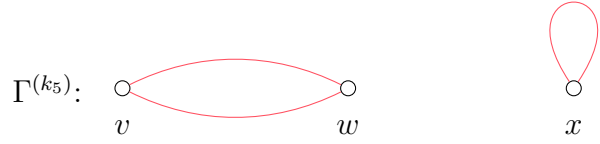
\threevertexgraph000101010{\Gamma^{(k_1)}}
\threevertexgraph111000000{\Gamma^{(k_2)}}
\threevertexgraph100000011{\Gamma^{(k_3)}}
\threevertexgraph010001100{\Gamma^{(k_4)}}
\threevertexgraph001110000{\Gamma^{(k_5)}}
The corresponding contributions are
\begin{align*}
    \f16D_{i,i}^3C(\Gamma^{(k_1)})\int_{\Gamma^{(k_1)}}f&=\f43D_{i,i}^3\int_{[0,1]^3} f(s_1,s_2)f(s_2,s_3)f(s_3,s_1) \Id s_1\Id s_2\Id s_3,\\
    \f16D_{i,i}^3C(\Gamma^{(k_2)})\int_{\Gamma^{(k_2)}}f&=\f16D_{i,i}^3\lb\int_{[0,1]}f(s,s)\Id s\rb^3,\\
    \f16D_{i,i}^3C(\Gamma^{(k_3)})\int_{\Gamma^{(k_3)}}f&=\f13D_{i,i}^3\int_{[0,1]^2} f(s_1,s_2)^2 \Id s_1\Id s_2\int_{[0,1]} f(s,s) \Id s,\\
    \f16D_{i,i}^3C(\Gamma^{(k_4)})\int_{\Gamma^{(k_4)}}f&=\f13D_{i,i}^3\int_{[0,1]^2} f(s_1,s_2)^2 \Id s_1\Id s_2\int_{[0,1]} f(s,s) \Id s,\\
    \f16D_{i,i}^3C(\Gamma^{(k_5)})\int_{\Gamma^{(k_5)}}f&=\f13D_{i,i}^3\int_{[0,1]^2} f(s_1,s_2)^2 \Id s_1\Id s_2\int_{[0,1]} f(s,s) \Id s.\\  
\end{align*}
The choice of $\alpha,\beta,\gamma$ which corresponds to $q_{2e_i}q_{2e_j}q_{2e_k}$ leads to only trivial sums except the $i$-th, $j$-th and $k$-th one. The only possible graph is
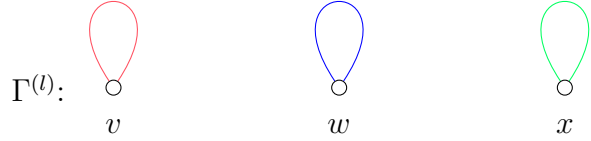
\threevertexgraph123000000{\Gamma^{(l)}}
The corresponding contribution is
$$
    \f16D_{i,i}D_{j,j}D_{k,k}C(\Gamma^{(l)}_{\t{red}})C(\Gamma^{(l)}_{\t{blue}})C(\Gamma^{(l)}_{\t{green}})\int_{\Gamma^{(l)}}f=\f16D_{i,i}D_{j,j}D_{k,k}\lb\int_{[0,1]}f(s,s)\Id s\rb^3
$$
The choice of $\alpha,\beta,\gamma$ which corresponds to $q_{e_i+e_j}q_{e_i+e_k}q_{e_j+e_k}$ leads to only trivial sums except the $i$-th, $j$-th and $k$-th one. The only possible graph is
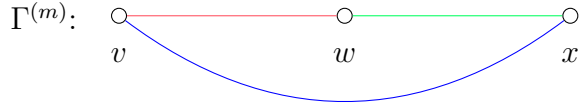
\threevertexgraph000102030{\Gamma^{(m)}}
The corresponding contribution is
$$
    \f16 8D_{i,j}D_{i,k}D_{j,k}C(\Gamma^{(m)}_{\t{red}})C(\Gamma^{(m)}_{\t{blue}})C(\Gamma^{(m)}_{\t{green}})\int_{\Gamma^{(m)}}f=\f43\int_{[0,1]^3} f(s_1,s_2)f(s_2,s_3)f(s_3,s_1) \Id s_1\Id s_2\Id s_3.
$$

\end{document}